\documentclass[11pt]{article}
\usepackage{style}
\usepackage{subcaption}
\usepackage{multirow}
\usepackage{array}
\usepackage{todonotes}
\newcommand{\rat}{\tau}

\newcommand*\samethanks[1][\value{footnote}]{\footnotemark[#1]}

\title{A comparative numerical study of graph-based splitting algorithms for linear subspaces}
\author{
    Francisco J. Arag\'on-Artacho\thanks{Department of Mathematics, University of Alicante, \textsc{Spain}. e-mail:~\href{mailto:francisco.aragon@ua.es}{francisco.aragon@ua.es}, \href{mailto:ruben.campoy@ua.es}{ruben.campoy@ua.es}, \href{mailto:cesar.lopez@ua.es}{cesar.lopez@ua.es}}
    \and Rub\'en Campoy\samethanks[1]
    \and Irene L\'opez-Larios
    \and C\'esar L\'opez-Pastor\samethanks[1]
}
\date{}

\begin{document}

\maketitle

\begin{abstract}
    In this note, we test the performance of six algorithms from the family of graph-based splitting methods [\emph{SIAM J. Optim.}, 34 (2024), pp. 1569--1594]
    specialized to normal cones of linear subspaces. To do this, we first implement some numerical experiments to determine the best relaxation parameter for each algorithm. Then, we compare the number of iterations each algorithm requires to reach a given stopping criterion, using the previously identified best relaxation parameter. The numerical results allow us to identify some relevant patterns and provide numerical evidence that may guide further theoretical analysis.
\end{abstract}

\section{Introduction}

Let $U_i\subseteq\mathbb{R}^p$ be linear subspaces for $i=1,\ldots,n$. We are interested in solving the following feasibility problem:
\begin{equation}\label{eq:problem}
\text{Find }x\in\bigcap_{i=1}^n U_i.
\end{equation}
For the case $n=2$, this problem can be solved numerically using the celebrated \emph{Douglas--Rachford (DR) algorithm} \cite{LM79}. Given an initial point $v^0\in\mathbb{R}^p$ and a \emph{relaxation parameter} $\theta\in{]0,2[}$, the algorithm iterates according to the recurrence
\begin{empheq}[left = \empheqlbrace]{align*}
    x_1^{k+1}&=P_{U_1}\big(v^k\big),\\
    x_2^{k+1}&=P_{U_2}\big(2x_1^{k+1}-v^k\big),\\
    v^{k+1}&=v^k+\theta\big(x_2^{k+1}-x_1^{k+1}\big),
\end{empheq}
where $P_U$ denotes the projection onto a subspace $U\subseteq\R^p$.
The generated sequences verify that $v^k\to v^*$, for some $v^*\in\mathbb{R}^p$, and $x_i^k\to x^*=P_{U_1}(v^*)=P_{U_2}(2x^*-v^*)\in U_1\cap U_2$, for $i=1,2$. Furthermore, it was shown in  \cite{cosineDR} that the limit point can be expressed as
\begin{equation}\label{eq:FPcharDR}
v^*=P_{U_1\cap U_2}(v^0)+P_{U_1^\perp\cap U_2^\perp}(v^0).
\end{equation}
Consequently, it follows that \(x^* = P_{U_1 \cap U_2}(v^0)\), and thus the Douglas--Rachford algorithm
solves not only the original feasibility problem, but also a best approximation problem.

The linear convergence of DR when applied to subspaces was
analyzed in \cite{cosineDR,relaxDR}.
The authors proved that the convergence rate depends on both the relaxation parameter \(\theta\) and the \emph{cosine of the
Friedrichs angle} between \(U_1\) and \(U_2\), given by
\[
c_F(U_1,U_2)
:= \max_{\|u_1\|,\|u_2\| \leq 1}
\left\{
\langle u_1, u_2 \rangle
\;\middle|\;
u_1 \in U_1 \cap (U_1 \cap U_2)^\perp,\;
u_2 \in U_2 \cap (U_1 \cap U_2)^\perp
\right\}.
\]
In particular, it was proved that the optimal choice of the relaxation parameter is \(\theta = 1\),
for which the smallest convergence rate, precisely \(c_F(U_1,U_2)\), is attained.

The natural extension of this algorithm to more than two subspaces fails to converge in general (see, e.g.,~\cite[\S~3.3]{DRconv}). A standard way to overcome this difficulty, originally proposed by Pierra~\cite{pierra},
is to reformulate the feasibility problem in a suitable product space.
More precisely, problem~\cref{eq:problem} can be equivalently written as
\[
\text{Find } \mathbf{x} \in \mathcal{U} \cap \Delta_n,
\]
where
\begin{equation}\label{eq:PS}
\begin{aligned}
\mathcal{U} &:= U_1 \times \cdots \times U_n \subseteq (\mathbb{R}^p)^n, \\
\Delta_n &:= \{(x,\ldots,x) \in (\mathbb{R}^p)^n \mid x \in \mathbb{R}^p\}.
\end{aligned}
\end{equation}
This reformulation allows the DR algorithm to be applied again,
at the expense of working in a higher-dimensional space.

In recent years, several generalizations of the Douglas--Rachford algorithm to an arbitrary
number of sets, without relying on a product space reformulation,
have been independently proposed (see, e.g., \cite{ryu20,malitsky2023resolvent,campoy,condat2023proximal}).
In the seminal work \cite{graph-drs}, these approaches have been unified under a single theoretical
framework through the introduction of the class of \emph{graph-based splitting methods}. This approach, which is explained in \cref{sec:graph} for the context of linear subspaces, provides a flexible and general methodology
for designing splitting algorithms.

The purpose of this work is to provide a numerical study of the performance of several instances
of graph-based DR methods, detailed in~\cref{sec:alg}, when applied to feasibility problems
over linear subspaces.
In particular, we investigate the influence of the relaxation parameter in \cref{sec:relax}
and compare the performance of different graph-induced algorithms in \cref{sec:compare}. Finally, some conclusions and open questions are discussed in \cref{sec:conclusion}.

\section{Graph framework for projection algorithms}\label{sec:graph}

The framework introduced in~\cite{graph-drs} is based on the choice of two oriented graphs \(G\) and~\(G'\), where
\(G = (\{1,\ldots,n\},\mathcal{E})\) is a connected graph that preserves the natural order
(that is, \((i,j)\in\mathcal{E}\) implies \(i<j\)),
and \(G'= (\{1,\ldots,n\},\mathcal{E}') \subseteq G\) is a connected subgraph maintaining all nodes.

Given the graph \(G\) and a node \(i\in\{1,\ldots,n\}\), and denoting by $|\cdot|$ the cardinality of a set, we define
\[
d_i^{\mathrm{in}} := |\{h \mid (h,i)\in\mathcal{E}\}|
\quad \text{and} \quad
d_i^{\mathrm{out}} := |\{h \mid (i,h)\in\mathcal{E}\}|,
\]
as the \emph{in-degree} and \emph{out-degree} of node \(i\), respectively, and we set
\[
d_i := d_i^{\mathrm{in}} + d_i^{\mathrm{out}}
\]
for the total \emph{degree} of node \(i\). On the other hand, associated with the graph \(G'\), we construct its \emph{Laplacian} matrix \(L \in \mathbb{R}^{n\times n}\) given by
\[
L_{i,j} =
\begin{cases}
d'_i & \text{if } i=j,\\
-1, & \text{if } (i,j)\in\mathcal{E}' \text{ or } (j,i)\in\mathcal{E}',\\
0, & \text{otherwise},
\end{cases}
\]
where \(d'_i\) denotes the total degree of node \(i\) in the subgraph \(G'\). 
The Laplacian can be factored as $L=ZZ^*$, where $Z$ is a full-rank $n\times (n-1)$ matrix.
Such decomposition can always be obtained by means of singular value decomposition (see~\cite[Proposition~2.16]{graph-fb}).

With this notation, given initial points $v_1^0,\ldots,v_{n-1}^0\in\mathbb{R}^p$ and relaxation parameter $\theta\in{]0,2[}$, the graph-based DR algorithm introduced in \cite{graph-drs} iterates according to
\begin{empheq}[left = \empheqlbrace]{align*}
    x_i^{k+1}&=P_{U_i}\left(\frac2{d_i}\sum_{(h,i)\in\mathcal{E}}x_h^{k+1}+\frac1{d_i}\sum_{j=1}^{n-1}Z_{i,j}v_j^k\right),\ \forall i=1,\ldots,n\\
    v_j^{k+1}&=v_j^k-\theta\sum_{i=1}^nZ_{i,j}x_i^{k+1},\ \forall j=1,\ldots,n-1.
\end{empheq}
For further insight into the role of the graphs $G$ and $G'$ in the algorithm, we refer the interested reader to \cite[Remark~3.8]{graph-fb}.

Analogously to the classical DR algorithm, the generated sequences satisfy that $v_j^k\to v_j^*$ for all $j=1,\ldots,n-1$ and $x_i^k\to x^*$ with
\[
x^*
=
\frac{2 d_i^{\mathrm{in}}}{d_i} P_{U_i}(x^*)
+ \frac{1}{d_i} \sum_{j=1}^{n-1} P_{U_i}\!\left(Z_{i,j} v_j^*\right)\in \bigcap_{i=1}^n U_i,
\qquad \text{for } i=1,\ldots,n,
\]
see~\cite[Theorem~3.2]{graph-drs}.
Although this equation does not provide explicit information about the limit points $x^*$ and $\mathbf{v^*}$,  a closed-form
expression that generalizes~\cref{eq:FPcharDR} has been recently obtained in \cite[Theorem~4.8]{graph-linear}.
More precisely, let \(\boldsymbol{\delta}\in\mathbb{R}^n\) denote the \emph{degree-balance} vector of the graph \(G\),
defined by \(\delta_i := d_i^{\mathrm{in}} - d_i^{\mathrm{out}}\) for \(i=1,\ldots,n\),
and let \(\boldsymbol{\alpha}=(\alpha_1,\ldots,\alpha_{n-1})\in\mathbb{R}^{n-1}\) be the unique solution to the system
\(
Z \boldsymbol{\alpha} = \boldsymbol{\delta}
\).
Consider also the subspace \(E\subseteq (\mathbb{R}^p)^{\,n-1}\) given by
\begin{equation*}\label{eq:defE}
E
:=\left\{\,(e_1,\ldots,e_{n-1})\in(\mathbb{R}^p)^{n-1}\ \middle|\
\textstyle\sum_{j=1}^{n-1} Z_{i,j}\,e_j \in U_i^\perp,\ \forall i=1,\ldots,n \right\}.
\end{equation*}
Then the limit points $x^*\in\R^p$ and $\mathbf v^*=(v_1^*,\ldots,v_{n-1}^*)\in(\R^p)^{n-1}$ admit the explicit representation
\begin{subequations}\label{eq:explicit-limits-graphDR}
\begin{align}
x^*
&=
\,
P_{\cap_{i=1}^n U_i}\!\left(\frac{1}{\|\boldsymbol{\alpha}\|^2}\,\sum_{j=1}^{n-1} \alpha_j v_j^0\right),
\label{eq:explicit-xstar-graphDR}
\\[0.3em]
\mathbf v^*
&=\,
\left(\alpha_1 x^*,\ldots,\alpha_{n-1} x^*\right)
\;+\;
P_E(\mathbf v^0),
\label{eq:explicit-vstar-graphDR}
\end{align}
\end{subequations}
where $\mathbf v^0:=(v_1^0,\ldots,v_{n-1}^0)$. 

\section{Selected algorithms and experimental setup}\label{sec:alg}
\label{sec:algos-stopping}

In our computational study, we shall consider the following six algorithms from the family of graph-based DR methods:

\begin{itemize}\itemsep=-2pt
\item \emph{Sequential} (given in \cite{graph-drs}).
\item \emph{Complete} (proposed in \cite{graph-drs} for $n=3$, generalized  in \cite{graph-fb} for any $n$, allowing forward terms).
\item \emph{Parallel down} (see \cite{graph-drs}).
\item \emph{Parallel up} (studied in \cite{campoy,condat2023proximal}).
\item \emph{Malitsky--Tam} (see \cite{malitsky2023resolvent}).
\item \emph{Generalized Ryu} (proposed in \cite{ryu20} for $n=3$  and generalized for arbitrary $n$ in~\cite{graph-drs}).
\end{itemize}


All these methods can be obtained as particular instances of the graph-based
framework described in \cref{sec:graph}, differing only in the choice of the underlying
graph pair \((G,G')\). The specific graph configuration for each algorithm is detailed in \cref{tab:graph_algorithms}.

\begin{table}[ht]
\centering{\footnotesize
\renewcommand{\arraystretch}{1.15}
\setlength{\tabcolsep}{10pt}
\begin{tabular}{lll}
\hline
\textbf{Algorithm} & \textbf{Graph \(G\)} & \textbf{Subgraph \(G'\)} \\
\hline

\multirow{2}{*}{Sequential}
& Sequential
& Sequential \\
& \(\mathcal{E}:=\{(i,i+1)\mid i=1,\ldots,n-1\}\)
& \(\mathcal{E}':=\{(i,i+1)\mid i=1,\ldots,n-1\}\) \\
\hline

\multirow{2}{*}{Complete}
& Complete
& Complete \\
& \(\mathcal{E}:=\{(i,j)\mid 1\le i<j\le n\}\)
& \(\mathcal{E}':=\{(i,j)\mid 1\le i<j\le n\}\) \\
\hline

\multirow{2}{*}{Parallel down}
& Parallel down
& Parallel down \\
& \(\mathcal{E}:=\{(i,n)\mid i=1,\ldots,n-1\}\)
& \(\mathcal{E}':=\{(i,n)\mid i=1,\ldots,n-1\}\) \\
\hline

\multirow{2}{*}{Parallel up}
& Parallel up
& Parallel up\\
& \(\mathcal{E}:=\{(1,i)\mid i=2,\ldots,n\}\)
& \(\mathcal{E}':=\{(1,i)\mid i=2,\ldots,n\}\) \\
\hline

\multirow{2}{*}{Malitsky--Tam}
& Ring
& Sequential \\
& \(\mathcal{E}:=\{(i,i+1)\mid i=1,\ldots,n-1\}\cup\{(1,n)\}\)
& \(\mathcal{E}':=\{(i,i+1)\mid i=1,\ldots,n-1\}\) \\
\hline

\multirow{2}{*}{Generalized Ryu}
& Complete
& Parallel down \\
& \(\mathcal{E}:=\{(i,j)\mid 1\le i<j\le n\}\)
& \(\mathcal{E}':=\{(i,n)\mid i=1,\ldots,n-1\}\) \\
\hline
\end{tabular}}
\caption{Graph-generated algorithms considered in the numerical experiments.}
\label{tab:graph_algorithms}
\end{table}

In all numerical experiments, the ambient space was fixed to \(\mathbb{R}^{50}\).
All subspaces \(U_i\subseteq\mathbb{R}^{50}\) and initial points were randomly generated,
and the same set of starting points was used to initialize all algorithms for each problem instance.

\begin{figure}[ht!]
\centering
\begin{subfigure}{0.41\textwidth}
  \centering
  \includegraphics[height=6.5cm]{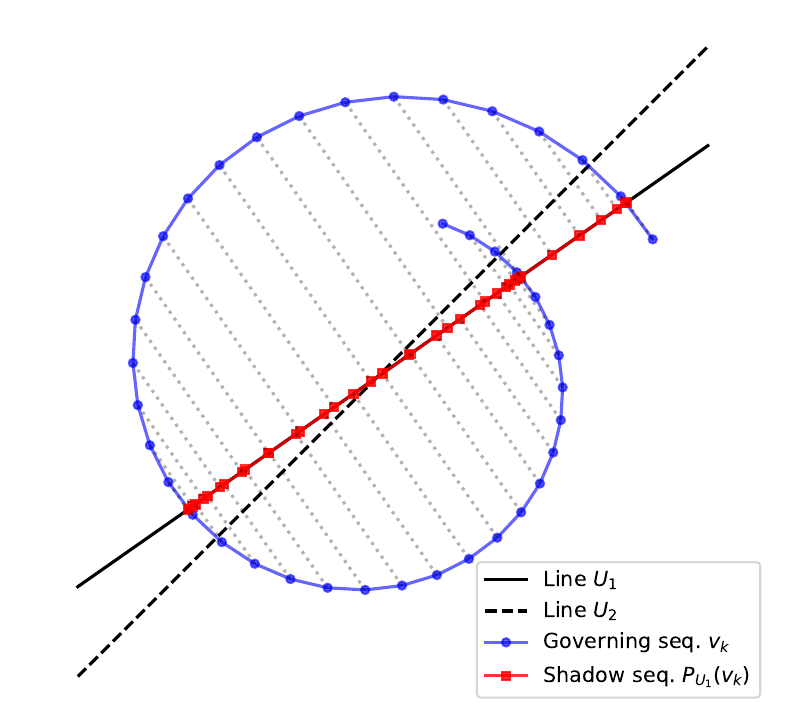}
  \caption{Trajectory of iterates}
\end{subfigure}\hfill
\begin{subfigure}{0.58\textwidth}
  \centering
  \includegraphics[height=6.5cm]{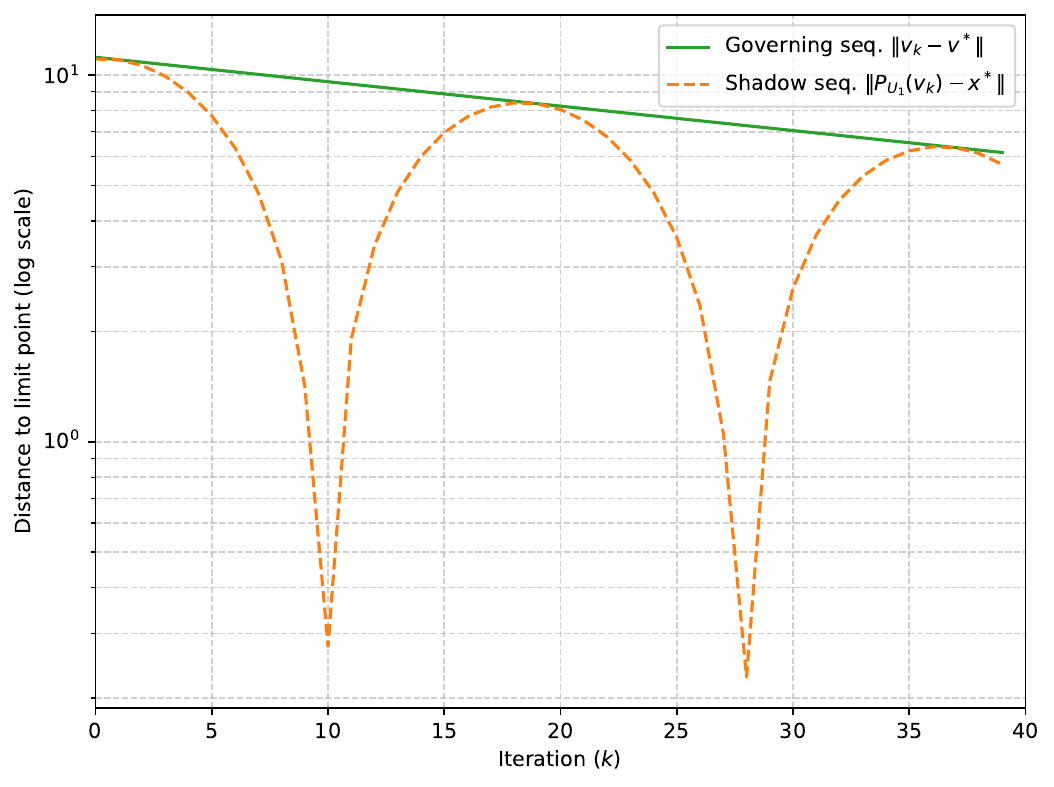}
  \caption{Distance to limit points}
\end{subfigure}
\caption{Douglas--Rachford algorithm for solving a feasibility problem involving two lines in \(\mathbb{R}^2\)}
\label{fig:dr-spiral}
\end{figure}

Although the shadow sequence is the one that converges to a solution $x^*$ of the problem,
in the numerical experiments we monitor the convergence of the governing sequence.
The reason for this choice is that the shadow iterates may happen to be close to the limit
even when the sequence has not yet converged. This phenomenon is already known for the classical Douglas--Rachford algorithm, where the governing iterates may exhibit a characteristic spiraling around the limit point, which results in a non-monotone behavior
of the distance between the shadow iterates and the solution
(see~\cref{fig:dr-spiral}).

Since the limit point can be explicitly computed
via~\cref{eq:explicit-limits-graphDR}, we adopted the following stopping criterion
in all numerical experiments:
\[
\|\mathbf v^k - \mathbf v^*\| < 10^{-6},
\]
where $\mathbf v^k:=(v_1^k,\ldots,v_{n-1}^k)$ is the governing sequence.
This choice provides a reliable measure of convergence and ensures a fair and consistent
comparison of the convergence behavior across different algorithms and parameter settings.

%

\section{Numerical study of the best relaxation parameter}\label{sec:relax}

For each $n \in\{ 3, 4, \ldots, 12\}$, we generated 20 random feasibility problems, each involving $n$ subspaces. For each problem instance, all algorithms were run from the same 10 randomly chosen starting points and the results were averaged. The algorithms were tested with different relaxation parameters $\theta$ selected from a grid $\{0.1, 0.2, \ldots, 1.9\}$.

\begin{figure}[htp!]
\centering
\includegraphics[width=\textwidth]{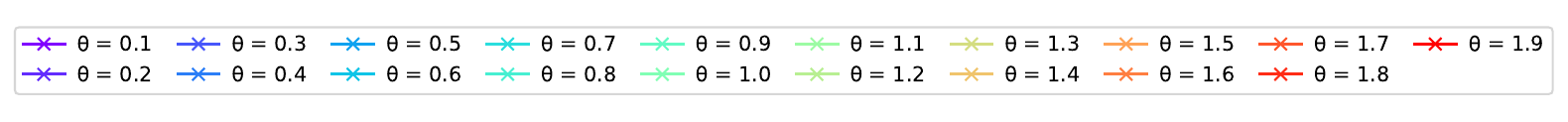}
\vspace{1mm}

\begin{subfigure}{0.49\textwidth}
  \centering
  \includegraphics[width=\linewidth]{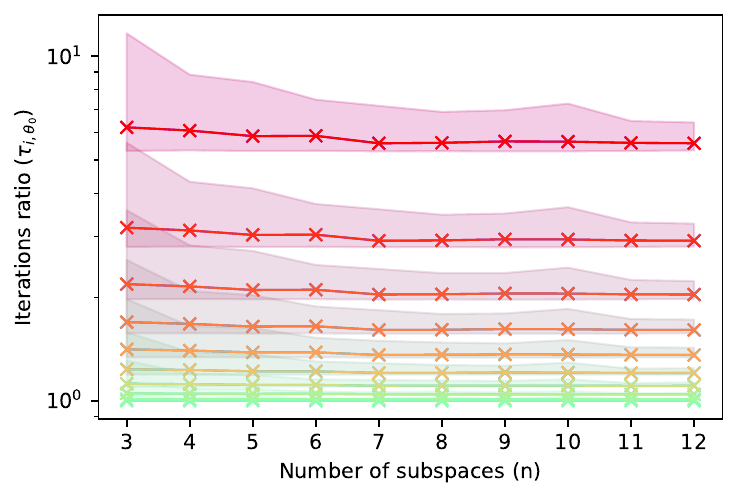}
  \caption{Sequential}
\end{subfigure}\hfill
\begin{subfigure}{0.49\textwidth}
  \centering
  \includegraphics[width=\linewidth]{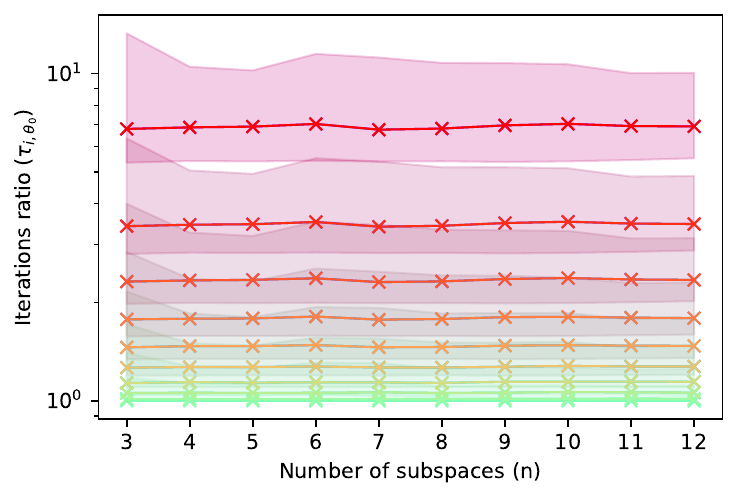}
  \caption{Complete}
\end{subfigure}
\vspace{5mm}

\begin{subfigure}{0.49\textwidth}
  \centering
  \includegraphics[width=\linewidth]{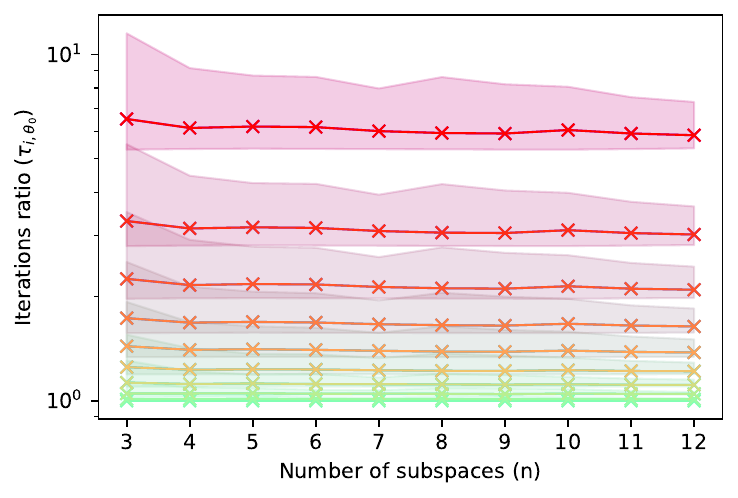}
  \caption{Parallel down}
\end{subfigure}\hfill
\begin{subfigure}{0.49\textwidth}
  \centering
  \includegraphics[width=\linewidth]{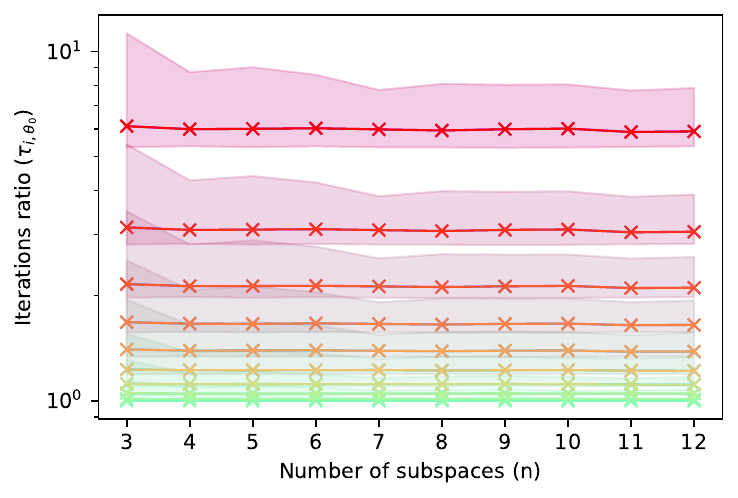}
  \caption{Parallel up}
\end{subfigure}
\vspace{5mm}

\begin{subfigure}{0.49\textwidth}
  \centering
  \includegraphics[width=\linewidth]{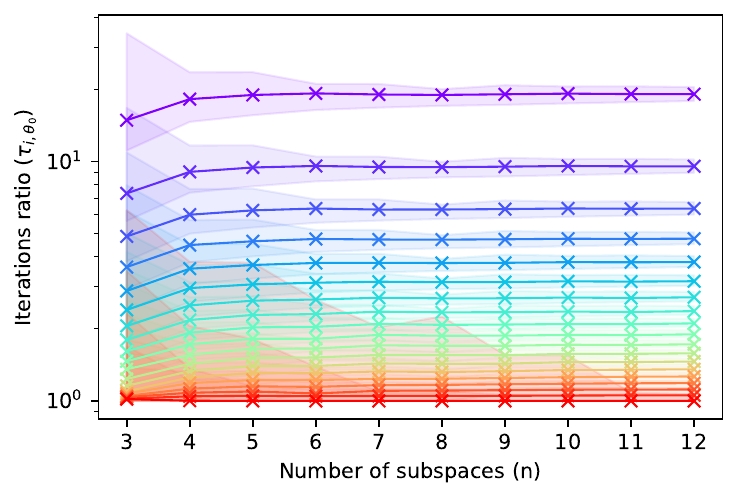}
  \caption{Generalized Ryu}
\end{subfigure}\hfill
\begin{subfigure}{0.49\textwidth}
  \centering
  \includegraphics[width=\linewidth]{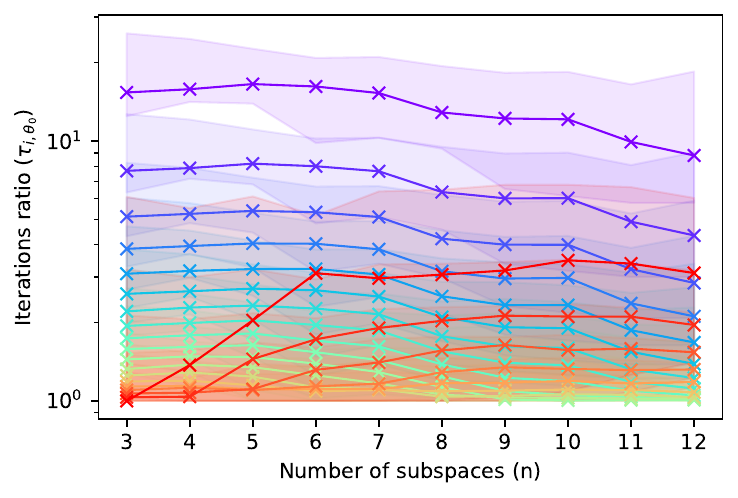}
  \caption{Malitsky--Tam}
\end{subfigure}
\caption{Results of the numerical experiment comparing the performance of each algorithm along different values of the relaxation parameter $\theta$}\label{thetas}
\end{figure}

Let $k_{i,\theta}$ denote the average number of iterations (among the 10 starting points) required by a given algorithm to solve the $i$-th feasibility problem using the relaxation parameter $\theta$. To assess the performance of a particular choice of parameter $\theta_0$, we compute the ratio
\[
\rat_{i,\theta_0} := \frac{k_{i,\theta_0}}{\min_{\theta \in \{0.1, \ldots, 1.9\}} k_{i,\theta}},
\]
which compares the number of iterations required by the algorithm using $\theta_0$ to the number required by the \emph{best-performing parameter} for that instance. The results for this test are shown in \cref{thetas}.  For each algorithm, a figure is plotted comparing the number of subspaces $n$ and the value of $\rat_{i,\theta_0}$. Varying colors depict different choices of $\theta_0$. Shaded areas represent the range of all possible ratios for all 20 instances, while the crosses connected by lines indicate the median of these ratios.

\paragraph{Discussion of results}
The first four algorithms depicted in \cref{thetas} behave in a very similar way, so we analyze them all together. Clearly, the best relaxation parameter
is $\theta=1$ (colored with light green).
There is another interesting fact. At first sight, one can observe the apparent absence of cool colors, i.e., results for values of $\theta\in{]0,1[}$. However, it turns out that these appear exactly behind the warm colors. Not only that, but the experiments show a symmetry between low and large values of $\theta$: the number of iterations required for $\theta$ and for $2-\theta$ coincides. In the plots, this is evidenced by the fact that, for example, the choices $0.1$ and $1.9$ completely overlap, as well as $0.2$ and $1.8$, and so on. These four algorithms yield identical results, suggesting that the inner graph structure explains these similarities. Notably, all four share the property $G=G'$.

Next, we observe some differences when examining the performance of the relaxation parameters in the Generalized Ryu algorithm. Indeed, now the best choice is $\theta=1.9$. 
On the other hand, the symmetry observed in the previous examples no longer applies here. In fact, a monotonic behavior is exhibited: the larger $\theta$, the better the performance.

Lastly, we have seen until now how the best choice for the relaxation parameter was independent of the number of subspaces. This is not the case for Malitsky--Tam: for small $n$, a large relaxation parameter is the best choice. However, as the number of subspaces increases, the value of the best-performing parameter reduces to reach the same value $\theta=1$ as in the four cases in which $G=G'$.

\begin{figure}[htp!]
\centering
\includegraphics[width=0.7\textwidth]{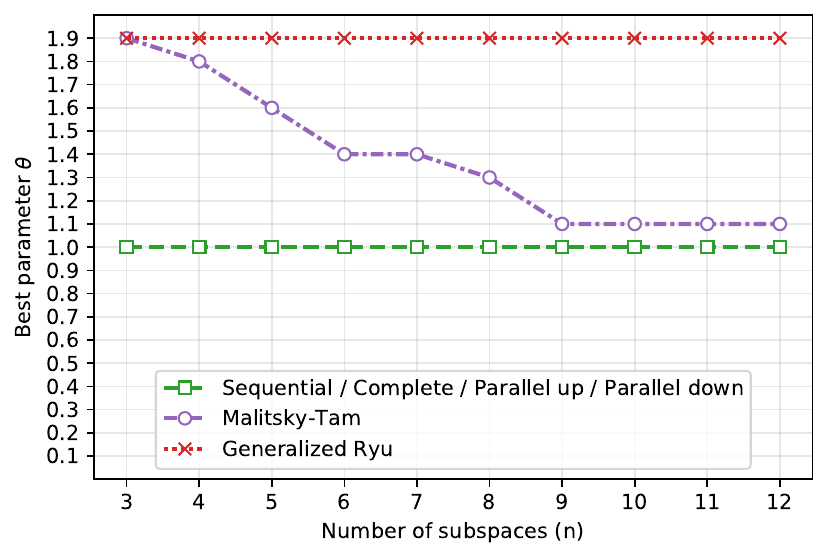}
\caption{Best relaxation parameter $\theta$ with respect to the number of subspaces}\label{best_theta}
\end{figure}

\paragraph{Overview} In \cref{best_theta}, the best-performing parameter for each number of subspaces is plotted for each of the six algorithms studied. The three categories described previously are clearly depicted: the four algorithms satisfying $G=G'$, whose best parameters are constantly equal to~$\theta=1$; Generalized Ryu, for which is always $\theta=1.9$; and the dependent-on-$n$ relaxation parameter of Malitsky--Tam. More precisely, we can observe that Malitsky--Tam's best parameter is close to (or even agrees with) the best parameter for Generalized Ryu when $n$ is small and, as $n$ increases, it eventually approaches to the best parameter for the other algorithms. A theoretical explanation of this behavior is also left as future work.

\section{Algorithmic performance}\label{sec:compare}

With the empirical results obtained in the previous section, we are now ready to fairly compare the performance of the six algorithms detailed in \cref{tab:graph_algorithms}, by setting for each of them the best-performing parameter according to~\cref{best_theta}. To this end, a similar experiment was conducted. For each $n\in\{3,4,\ldots,12\}$, representing the number of subspaces, we generated 100 random feasibility problems.
As in the previous experiment, we randomly generated 10 starting points on which the six algorithms were tested, and the results were averaged.

Since the convergence rate of the Douglas--Rachford algorithm is known to depend on the Friedrichs
angle~\cite{cosineDR,relaxDR}, it would be desirable to quantify this geometric measure also in the case of more than two subspaces.
A generalization of the Friedrichs angle to more than two subspaces, known as the \emph{Friedrichs number}, was introduced in \cite[Definition 3.2]{FriedrichsNumber}. This measure is explicitly determined by the Friedrichs angle between the associated subspaces under Pierra's product space reformulation given in~\cref{eq:PS}, see~\cite[Proposition~3.6(f)]{FriedrichsNumber}. For simplicity, in our numerical results we display Friedrichs angle under Pierra's reformulation. 
Although there is no theoretical guarantee that this quantity captures the geometry of
the original subspaces, the numerical results suggest that it provides a meaningful measure to do so.

\begin{figure}[htp!]

\centering
\includegraphics[width=0.6\textwidth]{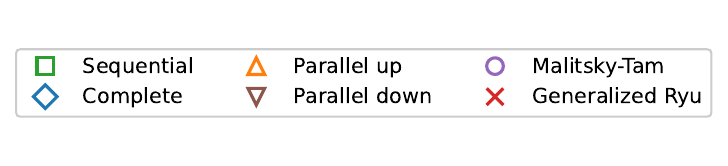}
\vspace{1mm}

  \begin{subfigure}{0.49\textwidth}
    \centering
    \includegraphics[width=\linewidth]{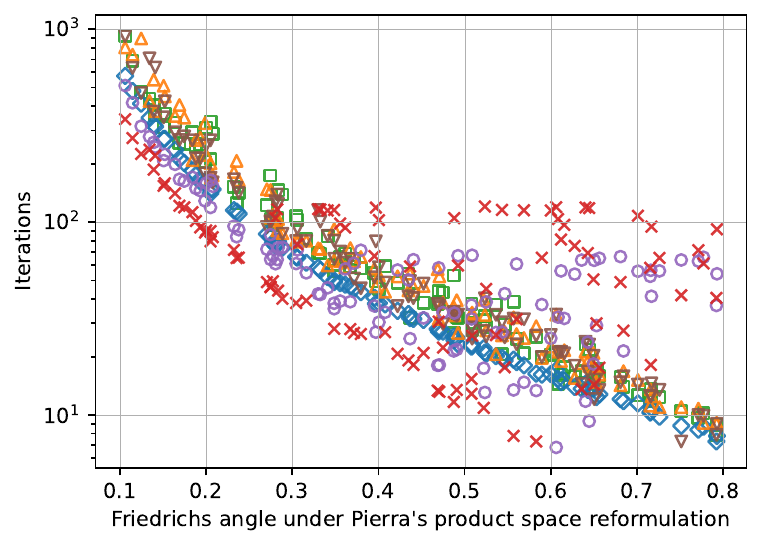}
    \caption{$n=3$}
  \end{subfigure}\hfill%
  \begin{subfigure}{0.49\textwidth}
    \centering
    \includegraphics[width=\linewidth]{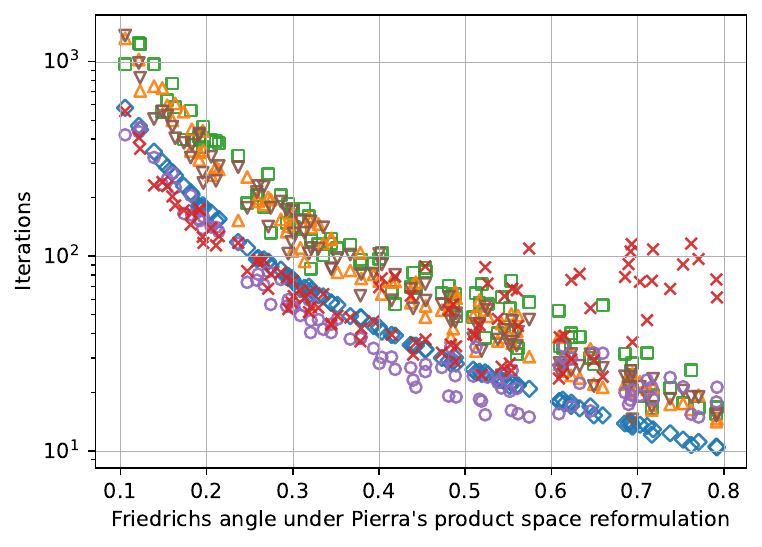}
    \caption{$n=5$}
  \end{subfigure}
\vspace{5mm}

  \begin{subfigure}{0.49\textwidth}
    \centering
    \includegraphics[width=\linewidth]{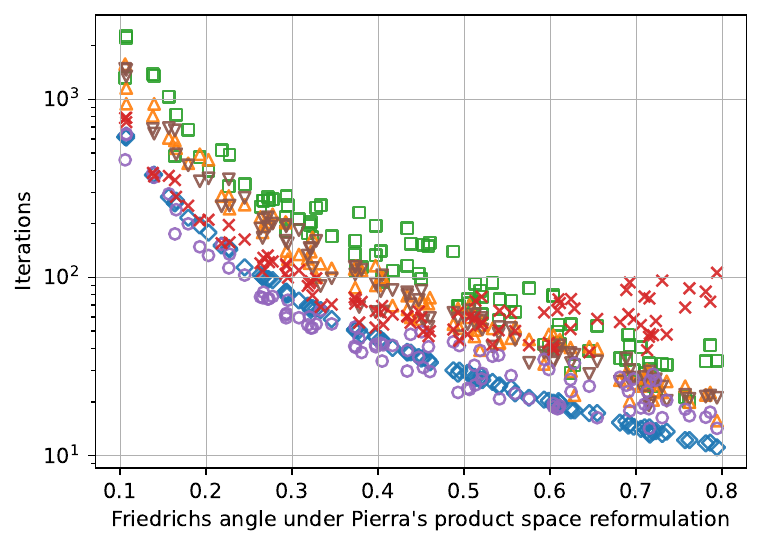}
    \caption{$n=7$}
  \end{subfigure}\hfill%
  \begin{subfigure}{0.49\textwidth}
    \centering
    \includegraphics[width=\linewidth]{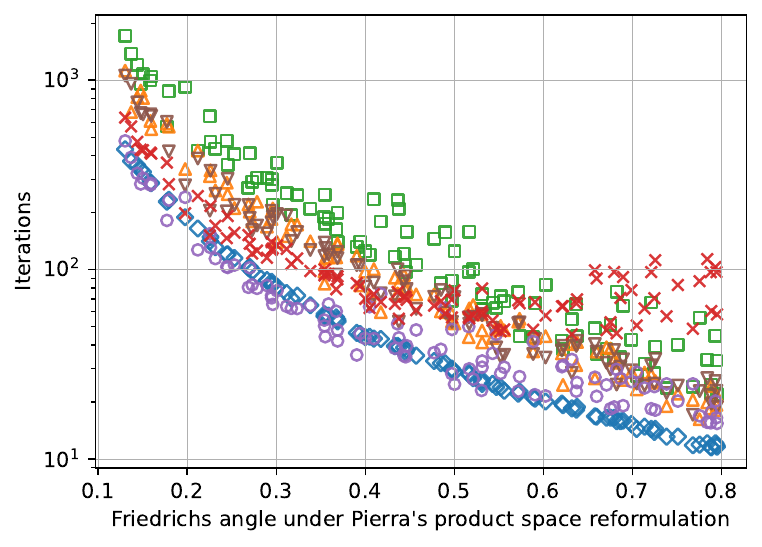}
    \caption{$n=8$}
  \end{subfigure}
\vspace{5mm}

  \begin{subfigure}{0.49\textwidth}
    \centering
    \includegraphics[width=\linewidth]{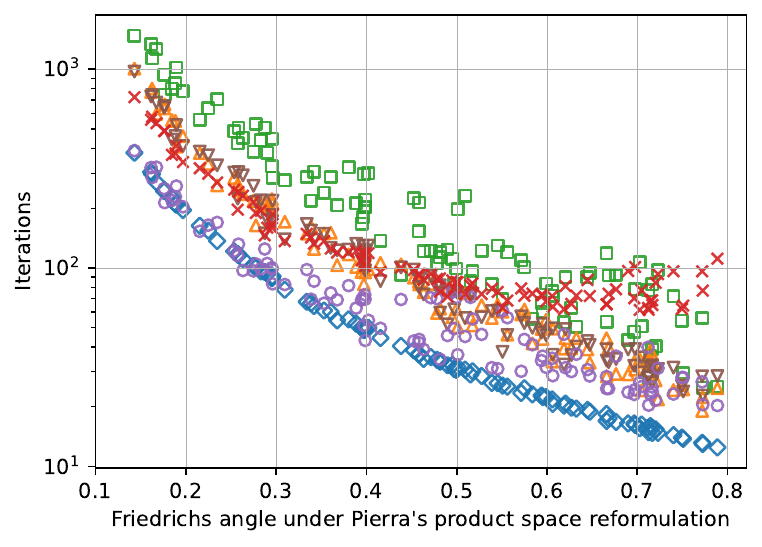}
    \caption{$n=10$}
  \end{subfigure}\hfill%
  \begin{subfigure}{0.49\textwidth}
    \centering
    \includegraphics[width=\linewidth]{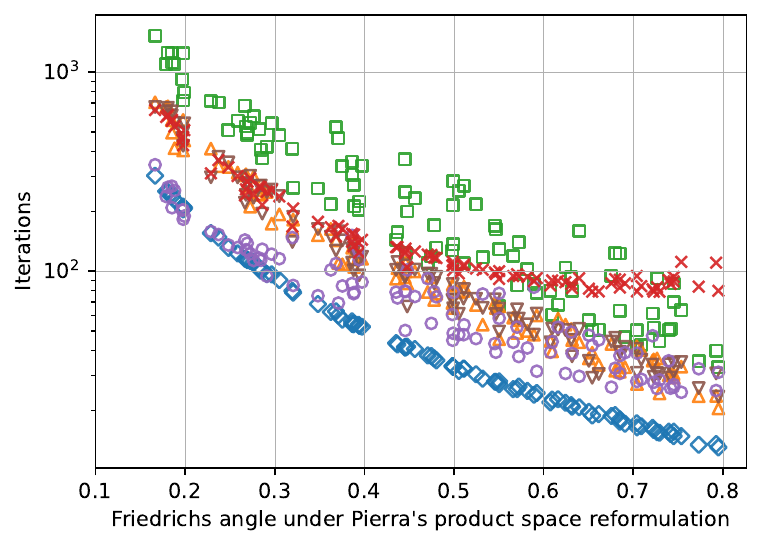}
    \caption{$n=12$}
  \end{subfigure}
\caption{Results of the numerical experiment comparing the performance of the algorithms with respect to the Friedrichs angle of the subspaces under Pierra's product space reformulation}\label{compare_angle}
\end{figure}

The results for this test, for selected values of $n\in\{3,5,7,8,10,12\}$, are shown in \cref{compare_angle}. Each plot illustrates the relationship between the Friedrichs angle (horizontal axis) and the mean number of iterations (vertical axis), with each marker representing the corresponding value for a given algorithm applied to a specific problem instance. Hence, there are 100 dots of the same shape and color for each algorithm in each plot. In our experiments, we computed Friedrichs angles from principal angles, see~\cite{relaxDR} for additional information.

\paragraph{Discussion of results}

We can observe in \cref{compare_angle} how the number of iterations is closely related to the Friedrichs angle. In particular, the relationship is very evident for the Complete algorithm.
For each number of subspaces $n$, the number of iterations appears to follow a curve that depends directly on the Friedrichs angle. The other algorithms seem to obey a similar curve pattern as the Complete, although the iterations are more spread and the relation is less evident. For instance, unlike the Complete, the Sequential algorithm suffers from this effect as $n$ grows, making less clear the relationship with the Friedrichs angle.

The Malitsky--Tam algorithm competes with the Complete in terms of performance. For small values of $n$ and small angles, the former was slightly faster, but it was slower for larger angles, showing a wider range and not adjusting well to a curve. On the other hand, for greater values of~$n$, the Complete algorithm was the best-performing method for nearly all problem instances.


The figure also illustrates how Parallel up and Parallel down are almost indistinguishable, using approximately the same number of iterations for each problem. This might be due to the graph configuration of both algorithms. Indeed, despite being defined by different algorithmic graphs, both parallel up and down graphs share the same graph topology: without taking the orientation into consideration, they are topologically isomorphic.

Lastly, an interesting behavior can be observed for the Generalized Ryu algorithm. Unlike the other algorithms, the typical curve fades out for larger angles.
This is particularly evident for
$n=3$, where, for angles greater than $0.5$, the number of iterations ranges from fewer than 10 to over 100. On the other hand, the curve becomes more apparent as $n$ increases. However, it flattens more than the other curves, meaning that the algorithm performed worse than the others.

\begin{figure}[htp!]
\centering
\includegraphics[width=0.85\textwidth]{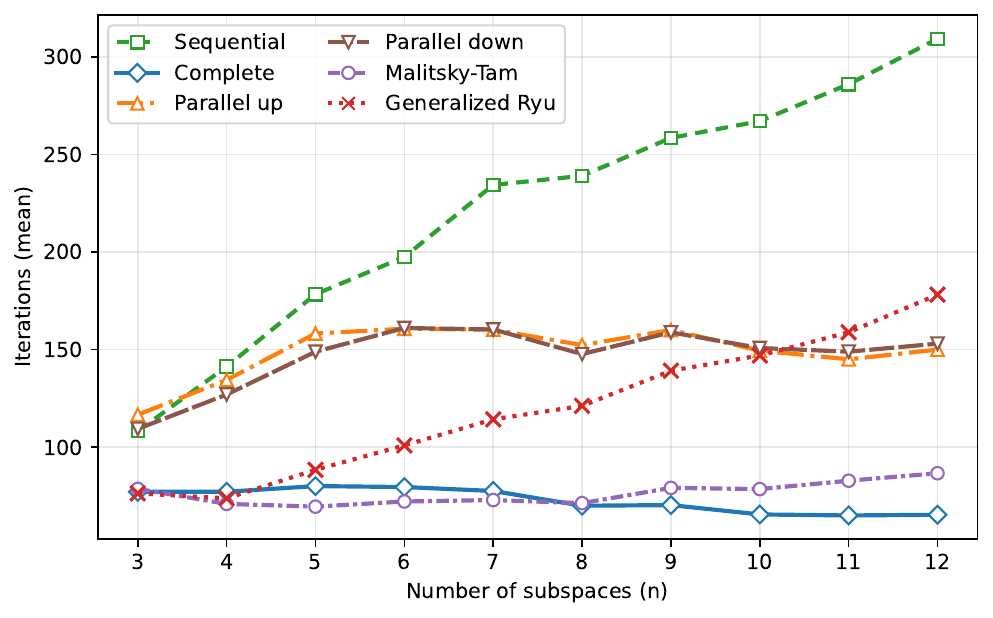}
\caption{Results of the numerical experiment comparing the performance of all the algorithms with respect to the number of subspaces}\label{compare_n}
\end{figure}

\paragraph{Overview} We summarize the information of the previous experiment in \cref{compare_n}. For each algorithm and each $n$, the mean of all the iterations is shown. This figure makes evident the most important aspects of our previous analysis. First, it shows that the Sequential algorithm was the slowest, with performance deteriorating as $n$ increases. The Parallel up and Parallel down algorithms follow, exhibiting nearly identical performance. Then comes the Generalized Ryu algorithm, whose performance worsens as $n$ increases. Finally, the nearly identical performance of the Complete and Malitsky--Tam algorithms is also illustrated in this plot. For smaller values of $n$, Malitsky--Tam required fewer iterations, whereas for $n\geq 8$, the Complete maintained the same number of iterations while Malitsky--Tam gradually required more to reach the stopping criterion.

\section{Conclusions and open questions}\label{sec:conclusion}
We have investigated the numerical performance of six characteristic algorithms from the recently introduced family of graph-based splitting methods~\cite{graph-drs}, applied to linear subspaces. We have numerically derived the best relaxation parameter for each algorithm and have compared them when the number of subspaces varies. The following questions have arisen:

\begin{itemize} \itemsep=-2pt
\item Is $\theta=1$ the optimal relaxation parameter when the graph $G$ is equal to the subgraph~$G'$? Why do $\theta$ and $2-\theta$ behave identically in this case?
\item Is there a general expression for the optimal relaxation parameter when $G\neq G'$?
\item Is it possible to obtain an expression for the rate of convergence in terms of the Friedrichs angle of the subspaces in the product space reformulation (or, equivalently, in terms of the Friedrichs number)?
\end{itemize} 
We believe that these questions deserve a theoretical explanation and plan to address them in future work.

\paragraph{Acknowledgments} We thank an anonymous reviewer for helpful comments. FJAA, RC and CLP were partially supported by Grant PID2022-136399NB-C21 funded by ERDF/EU. CLP was supported by Grant PREP2022-000118 funded by MICIU/AEI/10.13039/501100011033 and by ``ESF Investing in your future''.

\printbibliography

\end{document}